\documentclass[psamsfonts,onesided,10pt]{amsart}

\usepackage{amssymb,amsmath,amsfonts}
\usepackage[all,arc,cmtip]{xy}
\usepackage{tikz}
\usepackage{enumerate}
\usepackage{mathrsfs}
\usepackage[margin=1in]{geometry}
\usepackage{algorithm}
\usepackage{graphicx}
\usepackage{algpseudocode}
\usepackage{booktabs}
\usepackage{longtable}
\usepackage{wrapfig}

\theoremstyle{definition}

\theoremstyle{remark}

\makeatletter
\let\c@equation\c@thm
\makeatother
\numberwithin{equation}{section}

\title{Algorithms for Solving Optimization Problems Arising from Deep Neural Net Models: Smooth Problems}

\author{Vyacheslav Kungurtsev and Tomas Pevny}

\begin{document}
\maketitle

\begin{abstract}
Machine Learning models incorporating multiple layered learning networks have been seen to provide effective
models for various classification problems. The resulting optimization problem to solve for the optimal vector
minimizing the empirical risk is, however, highly nonlinear. This presents a challenge to application and development
of appropriate optimization algorithms for solving the problem. In this paper, we summarize the primary challenges
involved and present the case for a Newton-based method incorporating directions of negative curvature, including
promising numerical results on data arising from security anomally deetection. 
\end{abstract}

\section{Introduction}
In recent years there has been a notable increase in the popularity of models in machine learning incorporating multiple 
layers of classifiers, referred to as Deep Neural Nets (DNNs). 
Appled to text classification, perception and identification, and a myriad of other settings, deep learning
has, after a prolonged slow start, showed impressive efficacy~\cite{lecun2015deep}. Initially, the focus of optimization of these
problems was associated with the stochastic gradient descent (SGD) method~\cite{robbins1951stochastic}. This is because of two primary reasons.
First, SGD initially tended to outperform competing methods applied for large scale optimization problems arising from machine learning, 
as the historically common strongly convex models with strongly correlated
datasets permitted the low computational complexity per iteration of SGD to outweigh its relatively poor iteration complexity~\cite{bottou2016optimization}.  
Second, arising in the "big data" age, in which the size of datasets desired to be analyzed has grown exponentially, it is often impossible to calculate
and store an entire gradient vector, making stochastic and batch algorithms a necessity. 

However, it was quickly found that SGD as well as other first-order approaches were notoriously slow in solving these optimization 
problems~\cite{glorot2010understanding}. This can be attributed to the fact that the problems arising from DNNs are highly nonlinear. 
Originally, it was thought the primary challenge would be the presence of many local minima, necessitating techniques from global optimization,
which would prove challenging given their being subjected to the curse of dimensionality~\cite{zhigljavsky2007stochastic}. 
However, it was found that the more common
problem was the presence of saddle points, which were far more numerous and tended to drastically slow down the
performance first order algorithms~\cite{dauphin2014identifying}.

Second order methods incorporating a Newton direction can be used for large scale problems
arising in machine learning in a Hessian-free way~\cite{martens2012training}, however exact Hessian methods will encourage
algorithms to result in sequences moving towards saddle points, as they encourage rapid local convergence towards
any stationary point regardless of curvature. As such, Gauss-Newton type approximations to the
Hessian, or careful preconditioned steps are often used~\cite{pascanu2013revisiting, chapelle2011improved, vinyals2011krylov}, 
which encourage descent despite the presence of negetive curvature by essentially discarding the second derivative information
associated with the negative curvature. This is despite the fact that the actual direction along which negative curvature
is present gives a local second order approximation to the objective function that suggests a strong and reliable direction
of decrease. Direct use of negative curvature directions to calculate steps is limited~\cite{mizutani2008second}
despite well-developed iterative algorithms exploiting them in other areas of optimization~\cite{fasano2007iterative}.

\begin{figure}
\begin{center}
\includegraphics[scale=1.1]{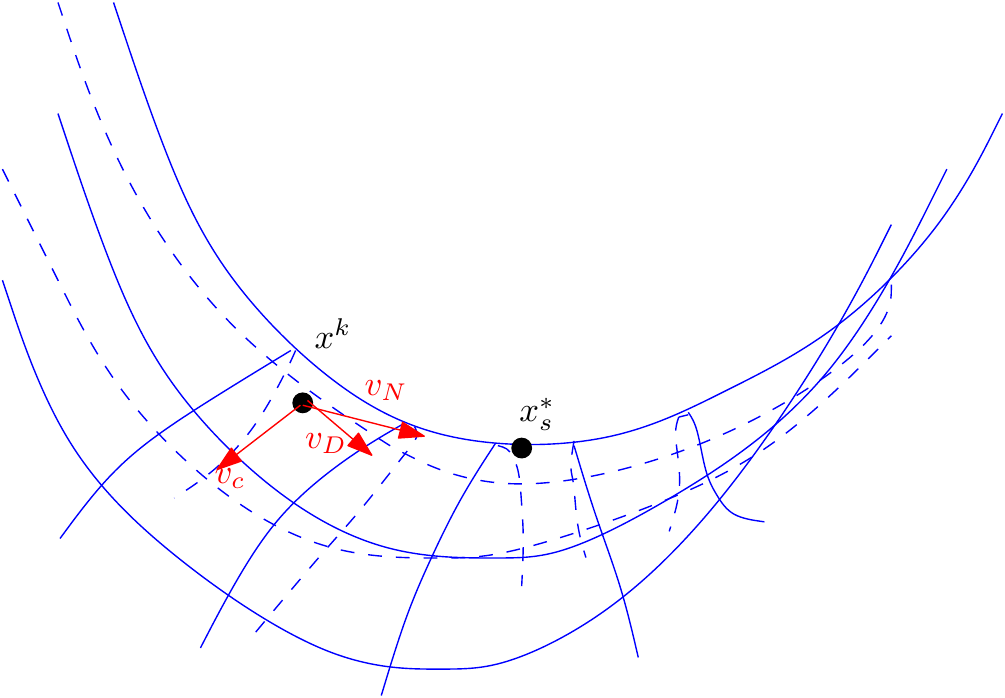}
\end{center}
\caption{\label{fig-saddle} Possible set of directions from a saddle point.}
\end{figure}

In order to motivate the use of direction of negative curvature, consider a generic saddle point in Figure~\ref{fig-saddle}. Here
$x^k$ is a current vector during an iterative optimization procedure, $x^*_s$ is the exact saddle point. A pure Newton direction will be
$v_N$, which will, if $x^k$ is sufficiently close to $x^*_s$, encourage rapid convergence of the iterates towards $x^*_s$ as the Newton
direction seeks first-order stationary points. A Gauss-Newton or regularized Newton direction will seek to create positive curvature in directions of negative curvature
and thus encourage weighted directions of descent, in this case labeled $v_D$. Here there is some movement away from the saddle point, but the
amount of curvature and hence level of scaling in the direction down the saddle is arbitrary and usually limited.
Finally, $v_c$ denotes the direction of negative curvature, which is clearly a movement away from the saddle point and thus, when the local
quadratic approximation is accurate, result in the largest decrease in the objection function.

\section{Lanczos Newton-Negative Curvature Method}
Consider minimizing a differentiable function $f(x):\mathbb{R}^N\to \mathbb{R}$, composed of a sum of component functions $f(x) = \sum_{i=1}^m f_i(x)$.

We consider incorporating a method incorporating Newton directions as well as directions of negative curvature using mini-batches. We use the
procedure presented in~\cite{lucidi1998curvilinear} to generate the direction of negative curvature and Newton curvature. The Lanczsos Algorithm
takes an initial vector and requires a mechanism of computing a Hessian vector product. Using finite differences, this can be approximated as,
\[
H(x^k)v \approx \frac{\nabla f(x^k+\epsilon v)-\nabla f(x^k)}{\epsilon}
\] 
using some scalar $\epsilon>0$. Alternatively, and more accurately, we can use a method akin to automatic differentiation~\cite{pearlmutter1994fast}
to compute the Hessian-vector product. 

The Lanczos algorithm takes an initial vector $g$ and the means of computing a Hessian vector with the Hessian evaluated at the current point
and returns matrices $V=(v_1\,v_2\,...\,v_q)$ and,
\[
T=\begin{pmatrix} \alpha_1 & \beta_2 & 0 & ... & 0 \\ \beta_2 & \alpha_2 & 0 & ... & ... \\ 0 & ... & ... & ... & ... \\ ... & ... & ... & ... & ... \\ 0 & ... & 0 & \beta_q & \alpha_q \end{pmatrix}
\] 
The vector $s = V y$ with $y$ such that $Ty = -V^T g$ can be considered an approximate Newton direction, in the sense that $s$ solves
$V^T (H(x^k)s+g)=0$. Alternatively, if we discard the columns of $V$ corresponding to indices $j$ for which $T_{jj}<0$ to get $\tilde V$, then the 
corresponding $\tilde s = \tilde V \tilde y$,  with $T \tilde y = -\tilde V^T g$, then we roughly get the direction corresponding to minimizing a quadratic 
incorporating only the directions of positive curvature. 

In addition, if $\mu$ is the minimum eigenvalue of $T$ and $w$ is the corresponding eigenvector, then $d = Vw$ is approximately a direction
of negative curvature for $H(x^k)$, in the sense that there exists an eigenvalue $\lambda$ of $H(x^k)$ such that $|\lambda-\mu|\le \beta_{q+1}$
where $\beta_{q+1}$ can be computed from the last steps in the Lanczos iteration.

Note that previous experiences with attempting to incorporate directions of negative curvature for problems arising in machine learning
have used the conjugate gradient (CG) procedure, which if prematurely terminated with the arising of a direction of negative curvature,
can then be used to calculate such a direction~\cite{mizutani2008second, martens2012training}. However, unlike the Lanczos procedure,
this method offers no guarantees as to the quality of this direction. 

We consider a generic Algorithm~\ref{alg:lan} which applies the Lanczos algorithm to a particular batch and takes a step
using the directions.

\begin{algorithm}[h]
	\caption{Batch Lanczos Subspace Descent Algorithm}
	\label{alg:lan}
	\begin{algorithmic}
		\State \textbf{Initialization:} $k=0$, $\mathbf{x}^0\in\mathbb{Re}^N$
		\While{$k\le k_{\text{max}}$}
		\State Choose index $j\in\{1,...,m\}$;
		\State{Compute $(s,\tilde s,d)=\text{LANCZOS}((H(x^k)(\cdot)),\nabla f_j(x^k))$;} 		
		\State{Compute the step $t$ as some combination of $(s,\tilde s, d)$;}
                \State{Compute a steplength $\alpha$}
		\State{Update the solution parameter vector $x^{k+1}=x^k+\alpha t$ and iteration counter $k \leftarrow k+1;$} 	
                \EndWhile
		\Return $\mathbf{x}^k$
	\end{algorithmic}
\end{algorithm}

Some details are in order:
\begin{enumerate}
\item The index $j$ can, potentially, be chosen in a number of ways, stochastically or deterministically. We found, however, in our
experiments that using round-robin or any other mini-batch determinstic cycling scheme works best, and random selection does not lead to a convergent procedure.
\item We note $H(x^k)(\cdot)$ as an argument for LANCZOS in the sense that we provide a way of computing 
the Hessian-vector product. We must also decide how many iterations the Lanczos procedure should perform, and after experiments investigating
the optimal amount, found that roughly $q=5$ iterations is best.
\item There are various forms of computing $t$, including $t=s+d$, $t=s+\frac{\|s\|d}{\|d\|}$, $t=\tilde s+d$, and so on. We found that, in general,
as long as both directions were used, it didn't appear that one form offered significant advantages over others, and thus we used the simplest $t=s+d$.
\item We use a line search procedure on the mini-batch function to calculate $\alpha$, i.e., we find $\alpha$ such that $f_j(x^k+\alpha t) < f_j(x^k)
+\eta\alpha t^T\nabla f_j(x^k)$, where $\eta>0$ is some small constant.  
\end{enumerate}

\section{Results}
We compared Algorithm~\ref{alg:lan} with standard stochastic gradient descent. We used SGD both with a constant, diminishing, and line search step size
and found that qualitatively the comparisons looked the same. For data we used private Skype security data arising from Cisco in Prague, Czech Republic. 
This is motivated by the fact that,
unlike many cases of classification model optimization where tight optimization tolerances and low risk measure objective values are weighed against the risk of 
overfitting~\cite{bousquet2008tradeoffs}, security cases requires low probability of false negatives above all else and thus it is desireable to push
objective values and tolerances as low as possible, with overfitting minimized purely by means of extensive cross-validation and careful selection 
of training and test sets. This setting creates a greater motivation for higher order methods. The security data is modeled with a mixture of Gaussian
layered in several layers.

\begin{figure}
\caption{\label{fig-convplot} Plot comparing the convergence of Algorithm~\ref{alg:lan} and SGD.}
\begin{center}
\begin{tabular}{c c} 
\includegraphics[scale=0.52]{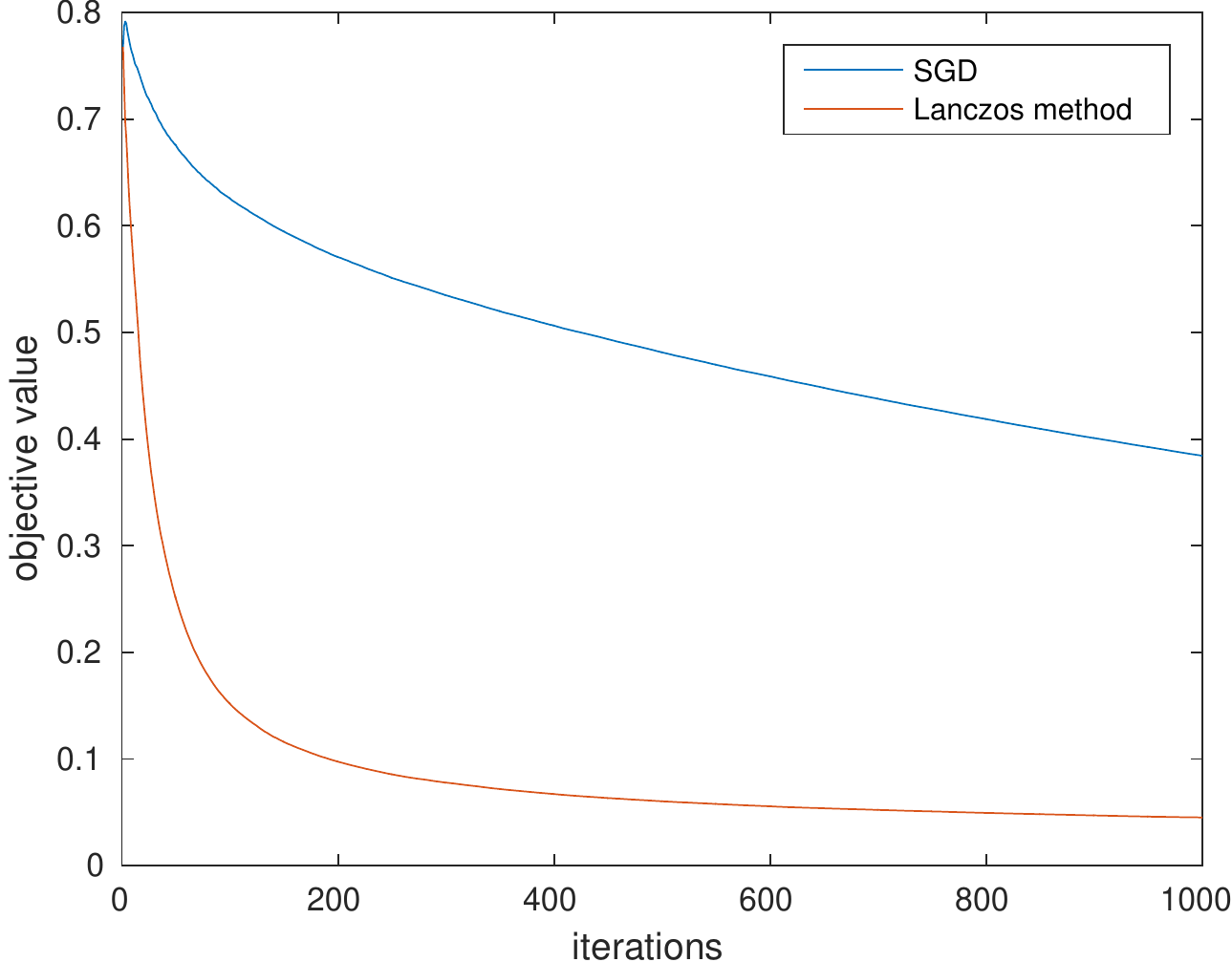}
&
\includegraphics[scale=0.31]{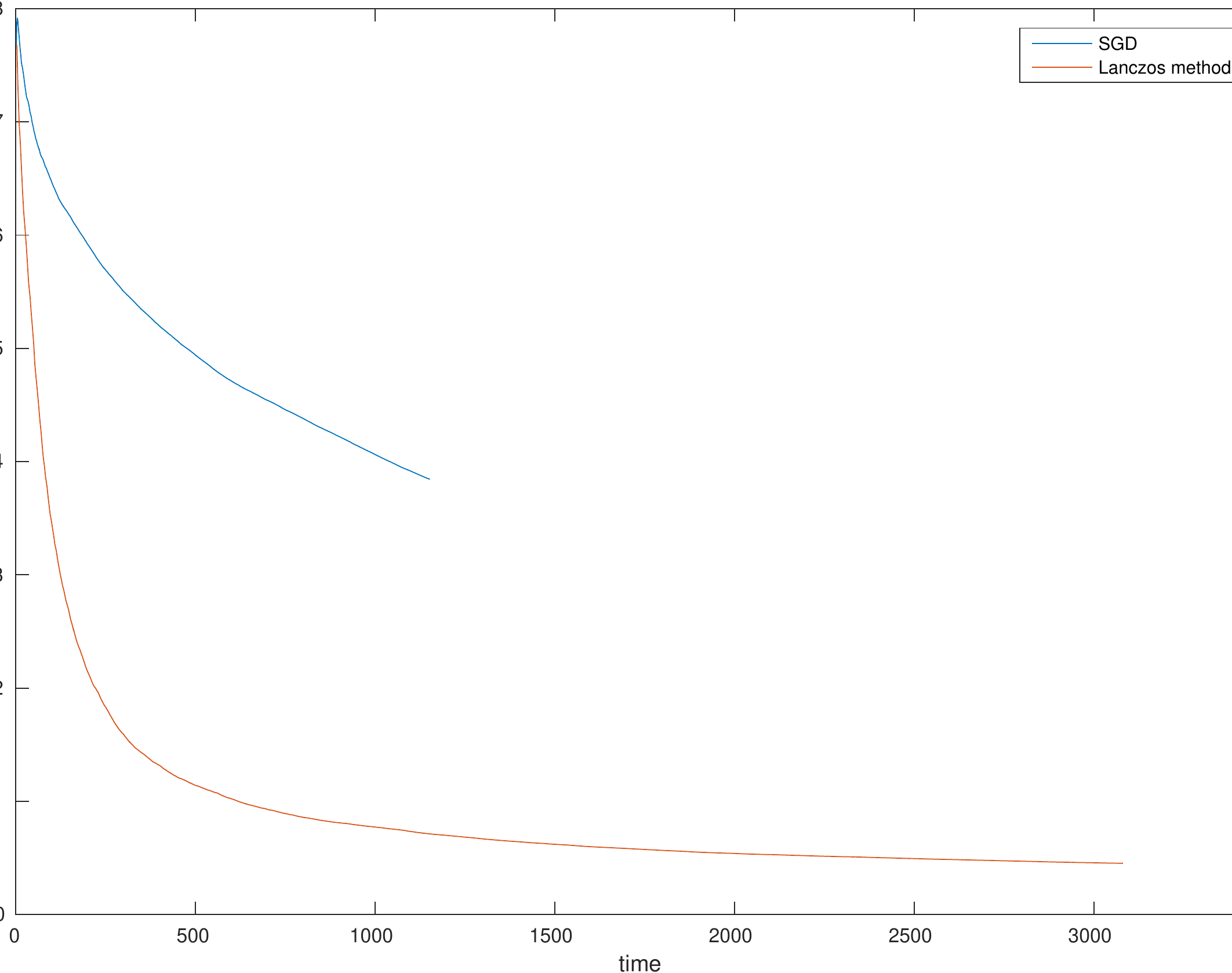}
\end{tabular}
\end{center}
\end{figure}

We ran 1000 iterations of each algorithm. The Lanczos algorithm, obviously, takes more time per iteration so we plot both the average objective
value over the iterations relative to the number of iterations as well as time. We see in Figure~\ref{fig-convplot} that Algorithm~\ref{alg:lan}
strongly outperforms SGD.

\section{Conclusion and Future Work}
This paper presents preliminary results suggesting that explicit calculation of directions of negative curvature, in order to take steps
down saddle points in the objective function map, are a potentially valuable algorithmic component for minimizing nonconvex machine learning
models arising from deep neural nets, by incorporating.

Aside from investigating the success of this method for a wider class of problems, there are two primary directions of future research
\begin{enumerate}
\item It has been found that second order information can be noisier than the gradient sample and still yield useful curvature information
for optimization. As such in stochastic Newton methods~\cite{byrd2011use} it is desireable to use a smaller batch size to compute the second order
information in order to decrease computation time per iteration. We will investigate how to carefully perform this.
\item We found deterministic cycling was essential to ensure convergence of the method, and random selection eventually resulted in objective
function ascent. We suspect this is due to random selection resulting in too high a chance of overfitting to particular mini-batches to the detriment
of other components of $f$. Other higher order algorithms incorporating mini-batches~\cite{sohl2014fast} tend to use historical information
from all components of $f$, regardless of when it was sampled, to continue to be used to update the vector. We will look to carefully implement
this, and see if it improves the performance of the stochastic method, although we note that the SFO algorithm in~\cite{sohl2014fast} also uses
round-robin index selection. 
\end{enumerate}

Incorporating these additional features will be helpful in advancing further the research into incorporation of negative curvature directions
for minimizing models arising from DNNs.

\bibliographystyle{plain}
\bibliography{refs}

\begin{thebibliography}{10}

\bibitem{bottou2016optimization}
L{\'e}on Bottou, Frank~E Curtis, and Jorge Nocedal.
\newblock Optimization methods for large-scale machine learning.
\newblock {\em arXiv preprint arXiv:1606.04838}, 2016.

\bibitem{bousquet2008tradeoffs}
Olivier Bousquet and L{\'e}on Bottou.
\newblock The tradeoffs of large scale learning.
\newblock In {\em Advances in neural information processing systems}, pages
  161--168, 2008.

\bibitem{byrd2011use}
Richard~H Byrd, Gillian~M Chin, Will Neveitt, and Jorge Nocedal.
\newblock On the use of stochastic hessian information in optimization methods
  for machine learning.
\newblock {\em SIAM Journal on Optimization}, 21(3):977--995, 2011.

\bibitem{chapelle2011improved}
Olivier Chapelle and Dumitru Erhan.
\newblock Improved preconditioner for hessian free optimization.
\newblock In {\em NIPS Workshop on Deep Learning and Unsupervised Feature
  Learning}, volume 201, 2011.

\bibitem{dauphin2014identifying}
Yann~N Dauphin, Razvan Pascanu, Caglar Gulcehre, Kyunghyun Cho, Surya Ganguli,
  and Yoshua Bengio.
\newblock Identifying and attacking the saddle point problem in
  high-dimensional non-convex optimization.
\newblock In {\em Advances in neural information processing systems}, pages
  2933--2941, 2014.

\bibitem{fasano2007iterative}
Giovanni Fasano and Massimo Roma.
\newblock Iterative computation of negative curvature directions in large scale
  optimization.
\newblock {\em Computational Optimization and Applications}, 38(1):81--104,
  2007.

\bibitem{glorot2010understanding}
Xavier Glorot and Yoshua Bengio.
\newblock Understanding the difficulty of training deep feedforward neural
  networks.
\newblock In {\em International conference on artificial intelligence and
  statistics}, pages 249--256, 2010.

\bibitem{lecun2015deep}
Yann LeCun, Yoshua Bengio, and Geoffrey Hinton.
\newblock Deep learning.
\newblock {\em Nature}, 521(7553):436--444, 2015.

\bibitem{lucidi1998curvilinear}
Stefano Lucidi, Francesco Rochetich, and Massimo Roma.
\newblock Curvilinear stabilization techniques for truncated newton methods in
  large scale unconstrained optimization.
\newblock {\em SIAM Journal on Optimization}, 8(4):916--939, 1998.

\bibitem{martens2012training}
James Martens and Ilya Sutskever.
\newblock Training deep and recurrent networks with hessian-free optimization.
\newblock In {\em Neural networks: Tricks of the trade}, pages 479--535.
  Springer, 2012.

\bibitem{mizutani2008second}
Eiji Mizutani and Stuart~E Dreyfus.
\newblock Second-order stagewise backpropagation for hessian-matrix analyses
  and investigation of negative curvature.
\newblock {\em Neural Networks}, 21(2):193--203, 2008.

\bibitem{pascanu2013revisiting}
Razvan Pascanu and Yoshua Bengio.
\newblock Revisiting natural gradient for deep networks.
\newblock {\em arXiv preprint arXiv:1301.3584}, 2013.

\bibitem{pearlmutter1994fast}
Barak~A Pearlmutter.
\newblock Fast exact multiplication by the hessian.
\newblock {\em Neural computation}, 6(1):147--160, 1994.

\bibitem{robbins1951stochastic}
Herbert Robbins and Sutton Monro.
\newblock A stochastic approximation method.
\newblock {\em The annals of mathematical statistics}, pages 400--407, 1951.

\bibitem{sohl2014fast}
Jascha Sohl-dickstein, Ben Poole, and Surya Ganguli.
\newblock Fast large-scale optimization by unifying stochastic gradient and
  quasi-newton methods.
\newblock In {\em Proceedings of the 31st International Conference on Machine
  Learning (ICML-14)}, pages 604--612, 2014.

\bibitem{vinyals2011krylov}
Oriol Vinyals and Daniel Povey.
\newblock Krylov subspace descent for deep learning.
\newblock {\em arXiv preprint arXiv:1111.4259}, 2011.

\bibitem{zhigljavsky2007stochastic}
Anatoly Zhigljavsky and Antanas {\v{Z}}ilinskas.
\newblock {\em Stochastic global optimization}, volume~9.
\newblock Springer Science \& Business Media, 2007.

\end{thebibliography}
\end{document}